\documentclass{svproc}

\usepackage{graphicx}
\usepackage{lscape}
\usepackage{amsmath}       % standard LaTeX graphics tool
\begin{document}

\title{Fibonacci and digit--by--digit computation;\\ An example of reverse engineering in\\ computational mathematics}

%\titlerunning{Computational Science in the 18th Century}

\author{Trond Steihaug} %\orcidID{0000-0003-1734-4638}
\institute{Department of Informatics, University of Bergen, Norway }
%\institute{ Department of Informatics, University of Bergen,
%Box 7803, N-5020 Bergen, Norway,\\
%\email{Trond.Steihaug@ii.uib.no}
%Department of Something, Somewhere\and The Third Person}

\titlerunning{Digit--by--digit computation}

\maketitle

%\pagestyle{myheadings}
%\thispagestyle{plain}
%\markboth{Historic development}{Nonlinear Equations}

\abstract{
The Fibonacci numbers are familiar to all of us.
They appear unexpectedly often in mathematics, so much there is an entire journal and a sequence of conferences dedicated to their study.  However, there is also another sequence of numbers associated with Fibonacci. In {\em The On-Line Encyclopedia of Integer Sequences}, a sequence of numbers  which is an approximation to the real root of the cubic polynomial. Fibonacci gave the first few numbers in the sequence  in the manuscript {\em Flos} from around 1215. Fibonacci stated an error in the last number and based on this error we try, in this paper to reconstruct the method used by Fibonacci. Fibonacci gave no indication on how he determined the numbers  and the problem of identifying possible methods was raised already the year after the first transcribed version of the manuscript was published in 1854. There are three possible methods available to Fibonacci to solve the cubic equation; two of the methods have been shown to give Fibonacci's result. In this paper we show that also the third method gives the same result and we argue that this is the most likely method.

}
\section{Introduction}
The Fibonacci numbers $1, 2, 3, 5, 8, 13, 21, 34, \ldots$ are well known \cite{Drozdyuk2010}. They are generated from the recurrence relation \(f_{i+2}=f_{i+1} + f_{i}, i=1,2,3\ldots\) where $f_{1}=1, f_2=2$. But there is another sequence of numbers associated with Fibonacci, namely $1,22,07,42,33,04,38, 30, 50 \ldots$. This sequence of numbers forms
an approximation to the real root of the cubic polynomial $x^3 + 2x^2 + 10x-20$ written in sexagesimal (base 60) notation. Fibonacci gave the numbers $1,22,7,42,33,4,40$ where the last sexagesimal digit is not correct.

Consider the cubic equation
\begin{equation}\label{eq:Fib}
x^3+2x^2+10x=20.
\end{equation}
From Cardano's formula, the only real solution can be found to be
\[\sqrt[3]{a+b} - \sqrt[3]{a-b} - \frac{2}{3}\]
where $a = 2\sqrt{3930}/9$ and $b=352/27$ \cite[p.21]{Posamentier2007}.

In {\em Flos} from around 1215 Fibonacci gives the numeric value in sexagesimal representation\footnote{In the Latin transcription by  Baldassarre Boncompagni \cite{Flos1854} from 1854}
{\em unum et minuta XXII et se
cunda VII et tertia XLII et quarta XXX\footnote{F.Woepcke \cite{Woepcke1854} pointed out that the fourth sexagesimal digit should be XXXIII and that the copyist has done similar errors in copying the manuscript.} et quinta IIII
et sexta XL} or in a classical sexagesimal notation
\begin{equation}\label{eq:Fib_number}
1^\circ\;22'\;7''\;42'''\;33^{IV}\;4^V 40^{VI}
\end{equation}
or
\begin{eqnarray*}
1+22/60+7/60^2+&42/60^3+33/60^4+4/60^5 +40/60^6 \approx&\\
&1.36880810785322371
\end{eqnarray*}
Fibonacci (or Leonardo of Pisa) was born sometime in the 1170s in Pisa, and died there sometime
after 1240. However, the existence of the manuscript {\em Flos} was not known until
Baldassarre Boncompagni published, a transcription in the original Latin language,  in 1854 \cite{Flos1854}. An algebraic version, using mathematical notation, was done by Angelo Genocchi in 1855 \cite{Genocchi1855}. Ettore Picutti gave a modern Italian translation with comments in 1983 \cite{Flos1983}.

The cubic equation (\ref{eq:Fib}) was known at the time of Fibonacci and is found in the work by Omar Khayyam (1048--1131)\footnote{A transcription and French translation in 1851 \cite[p.78]{Woepcke1851}, to English in 1931 \cite[p.110,p.114]{Kasir1931}
and German in 2012 \cite[p.155]{Linden2012} based on different manuscripts.}.

Woepcke \cite{Woepcke1854} in 1854 writes that the degree of accuracy of Fibonacci's
approximate value is  very remarkable, and  the knowledge of which method(s) used is of interest for the history of science. There are basically three types of methods Fibonacci with some certainty can have used. The first type of methods is secant type methods, the second type is Newton--Raphson type methods and the third type is based on digit--by--digit computation.

In a letter to Balthasar Boncompagni, Victor-Amédée Lebesgue (1791-1875) in 1855 argues that the method used by Fibonacci must have been a digit--by--digit method (the method of Vi\`{e}te)\footnote{Extract from the letter in Annanli di scienza matematiche e fisiche compilati da Barnaba Tortolini, Volume 6, 1855, p.155-160.}. Terquem in 1856 \cite{Terquem1856} states that Lebesgue advocates Vi\`{e}te's method.  Genocchi in 1855 determines one digit at the time in the sexagesimal (base 60) digit. However, it is done by determining the decimal representation of each sexagesimal digit which means up to two decimal digits for each sexagesimal digit.  Genocchi concludes that this is too labouriously  and suggests that Fibonacci must have used the golden rule of Cardano \cite[p.161-168]{Genocchi1855} which is a secant type method.

Hankel in 1874 \cite{Hankel1874} discusses a method for solving the cubic equation $Px=x^3+Q$ where $x_*^3$ is small compared to $P/Q$ which is the case here where $x_*=\sin(1)$\footnote{According to the arabian mathematicians M\={i}ram Čelebī (1475-1525) who described several methods for computing $\sin(1)$ known in Arabian world at the time of Fibonacci.}. The method is a digit--by--digit method and Hankel suggests that Fibonacci is using a similar method for (\ref{eq:Fib}) \cite[p.293]{Hankel1874}. This is further explored in Section \ref{sec:digit-by-digit}.

Cantor, in his lectures on the history of mathematics from 1892, suggests in addition to a digit--by--digit computation, that Cardano's golden method is a likely method to be used \cite[p.42--44]{Cantor1892}. Cardano's golden method is a secant type method.
%This is further explored in Section \ref{sec:interpolation}.
In the second edition from 1900 Cantor also discusses the approach taken by J.P Gram \cite[p.46-48]{Cantor1900} which is a Newton--Raphson type method.

Zeuthen \cite{Zeuthen1893} discusses possible methods Fibonacci may have used known at the time of Fibonacci and argues that Fibonacci may have used Newton-Raphson.

Gram \cite{Gram1893} shows four iterations with Newton-Raphson in sexagesimal arithmetic. Gram introduces truncation/rounding in the arithmetic and shows that this will give exactly Fibonacci's result \(1^\circ\;22'\;7''\;42'''\;33^{IV}\;4^V 40^{VI}\).

The last sexagesimal digit  is supposed to be $38^{VI}$.
The correct solution continues\footnote{Sequence A159990 in The On-Line Encyclopedia of Integer Sequences®(OEIS®)}
\[1^\circ\;22'\;7''\;42'''\;33^{IV}\;4^V\; 38^{VI}\; 30^{VII}\; 50^{VIII}\; 15^{IX}\; 43^{X}\; 13^{XI}.\]
Cassina in 1924 computes \(38^{VI}\; 30^{VII}\; 50^{VIII}\; 15^{IX}\) but incorrectly states the tenth and eleventh  sexagesimal digits \cite{Cassina1924}.

Al-Biruni  (973 – after 1050) states the equation $x^3=3x+1$ and gives the solution \cite[p.506]{Tropfke1980} and \cite[p.150]{Dattolini2021}.
\[1^\circ\;52'\;45''\;47'''\;13^{IV}\]
The correct solution is  \(1^\circ\;52'\;45''\;47'''\;12^{IV}\;43^V\;50^{VI}\). Al-Biruni result is an overestimate  similar to the overestimate made by Fibonacci in (\ref{eq:Fib}). It might be argued that since $43^V>30^V$ the fourth fractional digit is correctly rounded. The solution is $x_*=2\cos(\pi/9)$ \cite{Amir-Moez1994}.

David Eugene Smith \cite[p.471-472]{Smith1925} writes in 1925 and 1958 on Fibonacci's approximation:
\begin{quote}
No indication is given of the computational procedure to get this approximate solution.
How this result was obtained no one knows, but the
fact that numerical equations of this kind were being solved in
China at this time, and that intercourse with the East was
possible, leads to the belief that Fibonacci had learned of the
solution in his travels, had contributed what he could to the
theory, and had then given the result as it had come to him.
\end{quote}

David Eugene Smith is one of the founders of the fields of mathematics education and history of mathematics.

Glushkov in 1976 examines all works by Fibonacci and finds seven places where approximation methods are used \cite{Glushkov1976}. Glushkov argues that Fibonacci has used iterated linear interpolation to compute square and cube roots in {\em Liber abaci} Chapter 14 and suggests that Fibonacci made use of iterated linear interpolation also for  (\ref{eq:Fib}). He shows that starting with the two points 1 and 2 and using 18 iterations and truncation will give Fibonacci's numbers (\ref{eq:Fib_number}).  More recent interest in trying to explain the  error in Fibonacci's approximation  using two likely algorithms to solve the problem are in \cite{Brown2008,Maruszewski2009} from 2008 and 2009.

In Section \ref{sec:digit-by-digit} we first show the state of art of digit--by--digit computation in the mid--17th century by discussing Newton's annotation of Viet\`{e}'s book. We demonstrate the computation with sexagesimal arithmetic and rounding and show that Fibonacci's number (\ref{eq:Fib_number}) is a very likely result.

Zeuthen \cite{Zeuthen1893} argues that with the iterative methods available for Fibonacci, the Newton--Raphson method is very likely. In Section \ref{sec:NR} we discuss some published variations.

Most historians of mathematics argues that iterated linear interpolation is a plausible method used by Fibonacci to solve (\ref{eq:Fib}). This is discussed in Section \ref{sec:secant}.

In Section \ref{sec:geometric} we briefly state two geometric solution techniques known at the time  Fibonacci and show that the two sexagesimal digits may easily be found from the geometric construction.

\section{Digit--by--digit}\label{sec:digit-by-digit}

To illustrate the ideas behind the digit--by-digit computations we choose an example of depressed cubic polynomials (without the second order term) from a manuscript by Isaac Newton.
 In this unpublished note from 1664(?) reproduced in \cite[p. 63--71]{WhitesideI} Newton annotates Vi\`{e}te's
 {\em Opera Mathematica} from 1646 using the simplified notation in Oughtred's {\em Clavis Mathematic{\ae}} from 1648.
 Newton gives 7 examples  computing the root digit--by--digit. This unpublished note represents the 'state of the art'  in mid--17th century.
%\begin{itemize}
%\item For problem P.69  \(x^3-10x=13584\) Harriot--Commentary points out that Harriot in Petworth 9 is solving the same type of problem, but with different numerical terms.
%\end{itemize}
The  problems are
 on the form $x^3 + cx = d$ where $c$ and $d$ are positive real numbers. Define $p(x)= x^3 + cx - d$. Assume that the root is on the form
 $x = \alpha_2 10^2 +\alpha_1 10 + \alpha_0+ \alpha_{-1}/10$. A meta description of the algorithm is:
 \begin{itemize}
 %\item Step 1: Determine the number of digits in the root, say $x = \alpha_2 10^2 +\alpha_1 10 + \alpha_0$.
 \item Step 1: Determine the first digit $\alpha_2$: Choose the largest $0< \alpha_2\leq 9$ so that
 \[ p(\alpha_2 10^2) \leq 0
 \text{ and let } x_1 = \alpha_2 10^2.\]
 \item Step 2: Determine the second digit $\alpha_1$: Choose the largest $0\leq \alpha_1\leq 9$ so that
 \[p(x_1 + \alpha_1 10) \leq 0
 \text{ and let } x_2= x_1+\alpha_1 10.\]
 \item Step 3: Determine the digit   $\alpha_{0}$: Choose the largest $0\leq \alpha_{0}\leq 9$ so that
  \[p(x_2 + \alpha_0 ) \leq 0 \text{ and let } x_3= x_2+\alpha_0. \]
   \item Step 4: Determine the first digit in the fractional part  $\alpha_{-1}$: Choose the largest $0\leq \alpha_{-1}\leq 9$ so that
  \[p(x_3 + \frac{\alpha_{-1}}{10} ) \leq 0 \text{ and let } x_4= x_3+\frac{\alpha_{-1}}{10}. \]

 \end{itemize}
 The convergence of this technique follows from the observation that this is a bracketing process where the
 root will be in an interval on the form $[\cdot,\cdot)$ (the right
 end is open) and the monotonicity of $x^3 + cx$ in the interval.  The first interval
 will be $[\alpha_2 10^2,(\alpha_2+1) 10^2)$, then $[\alpha_2 10^2 +\alpha_1 10,\ \alpha_2 10^2 +(\alpha_1+1) 10)$.
 %and the final $[\alpha_2 10^2 +\alpha_1 10 + \alpha_0,\ \alpha_2 10^2 +\alpha_1 10 + \alpha_0+1)$.

 Stevin in 1585 suggests testing $0,1,2,3,\ldots$ to find the largest $\alpha_i$ \cite{Stevin1585}. Viet\`{e} in 1600 \cite{Viete1600} found an upper bound on $\alpha_i$ which would reduce the number of tests needed for find the largest $\alpha_i$. In the following few paragraphs, we show how to get an upper bound on the digits $\alpha_i$.

 Consider  \(p(x+h)=0\) for given $x>0$  where $p(x)<0$ and monotonically increasing. Then there exists unknown $h>0$
 \[p(x+h) = p(x) + h(3x^2+3xh + h^2 + c) = 0.\]
 Let $h\geq \underline{h}\geq 0$ be a lower bound on $h$, then an upper bound on $h$ will be
 \begin{equation} \label{Upper_bound}
 h \leq \frac{-p(x)}{3x^2+3x\underline{h}   + c}=\hat h,
 \end{equation}
 provided $c$ is not too negative, $-3x^2-3x\underline{h} <c$.
 Let $k$ be the number of digits in the integer part of the root. To determine digit number $j>1$, $\alpha _{k-j}$, consider
 \[x_j=\sum_{i=k-j}^{k-1} 10^{i} \alpha _{i} = x_{j-1} + 10^{k-j} \alpha _{k-j}.\]
 Define $h_j = 10^{k-j}$ then from (\ref{Upper_bound})
 \begin{equation}\label{Upper_bound_1}
 \alpha_{k-j} 10^{k-j} \leq \; \frac{-p(x_{j-1})}{3x_{j-1}^2+3x_{j-1}h_j  + c}
 \end{equation}
 and we have
 \begin{equation}\label{Upper_bound_2}
 \alpha_{k-j} \leq \Big \lfloor\; 10^{j-k}\quad \frac{-p(x_{j-1})}{3x_{j-1}^2+3x_{j-1}h_j  + c}\; \Big \rfloor.
 \end{equation}
 The  lower bound on $h$ to determine $ \alpha_{k-j}$ is in the case of  Vi\`{e}te $10^{k-j}$ while
 Holdred and Horner choose 0.

 In Table \ref{table:Newton_annotation} shows the actual computation of the digits in the solution. In the table the magnitude of the root, $k=3$, and the first digit, $2$, are known.
 \begin{table}[ht]
 \begin{center}
 \begin{tabular}{ |c|c |c | }
  \hline
              &\multicolumn{2}{|c|}{\(x^3+30x-14356197 \)}\\
  $k=3$       &\multicolumn{2}{|c|}{ \(x_*=243\)         }\\\hline
              \hline
              & $j=2$ &$j=3$\\ \hline
$x_{j-1}$           &    200&240     \\ \hline
$-p(x_{j-1})$       &6350197&524997  \\ \hline
$\underline{h}$    &10   &1    \\\hline
$3x_{j-1}^2+3x_{j-1}\underline{h}+c$&126030 &173550 \\\hline
$\hat h$      &50.4&3.03\\\hline
$\alpha_{k-j}$    &4 &3\\\hline
\end{tabular}
\end{center}
\caption{Equation \(x^3+30x-14356197=0 \). Number of digits in integral part of the positive root is $k=3$ and the first decimal digit  is 2, $j=1$ using Vi\`{e}te's bound $\underline{h}$.}
\label{table:Newton_annotation}
\end{table}
Consider finding root $x_*$ of \(x^3+30x-14356197 \). The number of digits in $x_*$, if $(x_*)^3 >>30x_*>0$, will be the number of digits in \(\sqrt[3]{14\;356\;197} \) which is 3 and the leading digit will be $\alpha_2=2$. To find the next digit of $x_*$ use
(\ref{Upper_bound}) with $x=200$ and $\underline{h}=10$. An upper bound of $\alpha_1\leq \lfloor 5.04\rfloor=5$. However, 5 is too large, $p(250) >0$, and the second digit is found to be 4,
$p(240) <0$, and $x_2=240$.

 Newton's transcripts of {Vi\`{e}te's solution of \(x^3+30x=14356197 \) is found in \cite[p.66]{WhitesideI} and reproduced in \cite[p.534]{Ypma1995}. The notebook (MS Add. 4000) with transcripts is available online\footnote{
https://cudl.lib.cam.ac.uk/view/MS-ADD-04000/1. Checked 01.04.2022.}.

Let $p(x) = x^3 + 2x^2 +10x-20$. Binomial expansion up to the third order was known at the time of Fibonacci \cite{Coolidge1949,Yadegari1980} so we may write
\begin{equation}\label{eq:Binomial_expansion}
p(x+h) = p(x)+h(3x^2+4x+10)+h^2(3x+2) + h^3
\end{equation}
The function $p(x)$ (for real $x$) is strictly monotone so if $p(x+h) \leq 0$ and for $x\geq 3/2$ and $0\leq \underline{h}\leq h$ we have an upper bound on $h$
\begin{equation}\label{eq:bound}
h \leq \frac{-p(x)}{3x^2+4x+10+\underline{h}(3x+2) + \underline{h}^2}
\end{equation}
This gives rise to four different divisors in (\ref{eq:bound})
\begin{enumerate}
\item Holdred and Horner--Ruffini: \(3x^2+4x+10\)
\item Viet\`{e}:  \(3x^2+4x+10+\underline{h}(3x+2)\)
\item Wallis: \(3x^2+4x+10+\underline{h}(3x+2) +\underline{h}^2\)
\item Hankel: A constant
\end{enumerate}
Note that the first divisor is $p'(x)$, the second is $p'(x)+\frac{\underline h}{2} p''(x)$ and the third is  $p'(x)+\frac{\underline h}{2} p''(x) + \frac{\underline h}{3!} p'''(x)$.

Fibonacci knew that a solution of (\ref{eq:Fib}) was between 1 and 2 so the first sexagesimal digit would be $1 (= 1^\circ)$. To determine the next sexagesimal digit, the first in the fractional part,  $\alpha_1$, using the first divisor in the list:
\[\alpha_1 \leq \left \lfloor 60\frac{-p(1)}{3\cdot 1^2+4\cdot 1+10}\right \rfloor = 24\]
using (\ref{eq:bound}) for $h=\alpha_1/60$.
However, the largest $\alpha_1$ so that $p(1+\alpha_1/60)\leq 0 $ is 22 so the second sexagesimal digit is $22 (= 22')$.
For the next digit $x_2=1+22/60$ and $h= \alpha_2/60^2$
\[ \alpha_2 \leq \left \lfloor 60^2\frac{-p(x_2)}{p'(x_2)} \right \rfloor = 7\]
and the approximation is $1^\circ 22'\;7''$. Assuming that we have computed $x=1^\circ\;22'\;7''\;42'''\;33^{IV}\;4^V$ and want to compute $\alpha_6$
(let $x_5=1+22/60+7/60^2+42/60^3+33/60^4+ 4/60^5$)
\[a_6 \leq \left \lfloor 60^6\frac{-p(x_5)}{p'(x_5)} \right \rfloor = 38.\]
If we do the computation by hand and using sexagesimal digits a possible rounding will be
\[-p(x_5) = 13^V32^{VI} \leq 14^V,\]
\text{ and }
\[p'(x_5)= 21^\circ\;5'\;46''\;6'''\;5^{IV}\;42^V 28^{VI}\geq 21^\circ.\]
Then
\[ \alpha_6 \leq 60^6\frac{14^V}{21^\circ}  = 40.\]
This may explain why the last digit in Fibonacci's approximation is 40. One might conjecture that Fibonacci just avoided the effort to calculate the value of $p(x_5+k/60^6))$ for $k = 38$ or $39$.

A change of variable $x+\frac{2}{3}$ in (\ref{eq:Fib}) gives the depressed cubic equation
\begin{equation}\label{eq:Fib_depressed}
x^3+\frac{26}{3}x=\frac{704}{27}
\end{equation}
which is considered by Hankel \cite{Hankel1874} and Gram \cite{Gram1893}.
Hankel's \cite{Hankel1874} digit-by-digit computation is probably best explained by defining $p(x) = x^3 + cx -d$, the form of cubic polynomials Newton is using.  Again, if $x>0$ and $h>0$ then
\[p(x+h) = p(x) + h(3x^2+3xh + h^2 + c)> p(x)+h(3x^2+c)  \]
and if the root is small, we can use the divisor $c$. Since the value of $x$ in the digit-by-digit computation is increasing and we know the first digit (hence $x_1$ is known) the constant divisor can be $3x_1^2+c$ and
\[  h \leq \frac{-p(x)}{3x_1^2+c}\text{ for } x\geq x_1.\]
For the depressed cubic (\ref{eq:Fib_depressed}) let $p_d(x) = x^3+\frac{26}{3}x-704/27$, $x_1=2^\circ$ and
$3x_1^2+8^\circ\;20' =20^\circ\;20'> 20^\circ$. The constant used by Vetter to approximate the derivative in Newton-Raphson \cite{Vetter1928} is $20^\circ$.
This simplifies the computation and assuming the first 5 sexagesimal digits are computed,
\[x_5 = 2^\circ\;2'\;7''\;42'''\;33^{IV}\;4^V\]
then
\[a_6 \leq \left \lfloor 60^6\frac{-p_d(x_5)}{20} \right \rfloor = 40.\]
Again, the approach taken by Hankel to use a constant divisor can explain that the 6th sexagesimal digit is 40 (and not 38).

   One might also conjecture that Fibonacci  calculated only $p(x_5+k/60^6)$ for $k = 36$ and $k = 40$, the less complex cases concerning sexagesimal fractional arithmetic\footnote{Reinhard Zumkeller, May 01 2009 in
the web-page of \em{The On-Line Encyclopedia of Integer Sequences}. Sequence A159990}?
$p(x_k+36/60^6) \approx -1.113673\; 10^{-7}$ and $p(x_5 + 40/60^6) \approx +6.719323\;10^{-10}$.
The latter result looks precise enough and can explain and justify Leonardo's 'rounding'.

Cassina \cite{Cassina1924} uses Horner-Ruffino's scheme and digit--by--digit to compute 14 decimal places. However, Cassina incorrectly states  the sexagesimal digits \(42^{X}\; 45^{XI}\) \cite{Cassina1924}.
The historian of mathematics Roland Calinger conjectures that Horner's method (a digit--by--digit method) was used \cite[p.369]{Calinger1999} by Fibonacci.

\subsection{Other sources of error in the last digit}
When computation became more documented in algebra books like in the 17th century \cite{Viete1600,Harriot1631,Oughtred1647}, the new value of the polynomial is shown to be computed iteratively
\(p(x+h)=p(x)+(p(x+h)-p(x))\). In general, the denominator in (\ref{Upper_bound_2}) and the difference $p(x+h)-p(x)$ and  will have common terms saving some arithmetic operations \cite{Steihaug2021}. Limited accuracy in $p(x)$ will thus propagate to the next $p(x+h)$. We know from Newton's manuscripts that only the number of digits needed are computed \cite{Ypma1995} and the most common way to terminate arithmetic operations is by truncating. The value of $p(x)$ must be reevaluated close to the end of the iteration sequence.

\section{Newton--Raphson Approaches}\label{sec:NR}

J{\o}rgen Pedersen Gram \cite[p.18-28]{Gram1893} shows that Fibonacci's result can be obtained by
Newton--Raphson.  Gram considers the depressed cubic (\ref{eq:Fib_depressed}) in sexagesimal notation
\[x^3+8^\circ 40'x = 26^\circ 4'\; 26{''}\;40{'''}.\]
Starting with $x_0=2^\circ$ and truncate the correction $2'$, the first iterate is $x_1=2^\circ 2'$. Truncating the next correction $8{''}$ and the second iteration is $x_2= 2^\circ 2'\; 8{''}$. Increasing the accuracy in the computation and computing the third correction $-17{'''}\;27^{IV}\;35^{V}$ and $x_3= 2^\circ2'\;7{''}\;42{'''}\;32^{IV}25^{V}$. The denominator (the derivative) can be replaced by $21^\circ 5'$ and the correction is $39^V\;40^{VI}$. The final approximation is then
\begin{eqnarray*}\nonumber
x_4 - 40' &=& 2^\circ\; 2'\;7''\;42'''\;33^{IV}\;4^V\; 40^{VI} - 40'\\ \nonumber
             &=& 1^\circ\;22'\;7''\;42'''\;33^{IV}\;4^V\; 40^{VI}
             \end{eqnarray*}
which is Fibonacci's approximation.

Gram shows that it is possible to get Fibonacci's result by approximating function and derivative values and truncate the computation.

Vetter \cite{Vetter1928} in 1928 replaces the denominator in the Newton--Raphson approach by the constant 20. However, Vetter advocates that the number of sexagesimal digits should not increase (in the approximate root) with more than one digit in each iteration. This will reduce the computational cost of doing calculation by hand. However, more than one sexagesimal digit can change from one iteration to the next. Figure \ref{fig:Vetter} shows the sexagesimal digits in the approximation in nine iterations
\begin{figure}[ht]
\begin{center}\includegraphics[scale=.6]{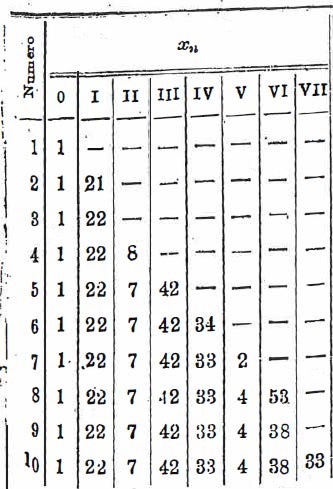}\end{center}
\caption{Sexagesimal computation of the approximation
\(1^\circ\;22'\;7''\;42'''\;33^{IV}\;4^V\; 38^{VI}\; 33^{VII}\) of the real root of $x^3+2x^2+10x-20=0$ in \cite{Vetter1928} using Newton-Raphson with a limit of at most one new digit in each iteration.}
\label{fig:Vetter}
\end{figure}
\section{Iterated linear interpolation}\label{sec:secant}
Fibonacci was confident in his problem-solving
abilities, and  Chapter 13 of Liber abaci, titled {\em Here Begins
chapter Thirteen on the Method Elchataym and How with It Nearly All Problems of
Mathematics Are Solved.} In this chapter Fibonacci explains the method {\em regulis elchatayn} (after al-kha\d{t}\={a}'ayn) that Fibonacci had learned from Arab sources \cite{Schwartz2004,Roberts2020}. The Arabic
al-khata’ayn is literally, “the two errors”, which is translated as the method of Double False
Position or Regula Falsi. A comprehensive history of the secant method/regula falsi/iterated linear interpolation is found in \cite{Papakonstantinou2013}. If $a$ and $b$ bracket a root of a concave  polynomial $p(x)$ (as is the case of $p(x) = x^3+2x^2+10x-20$ for $a=1(=x_{-1})$ and $b=2(=x_0)$) the new value, $x_1$, based on linear interpolation will be an underestimate of the root. This is repeated,  $x_1$ and $x_0$ bracketing the root. Fractional sexagesimal computation using a fixed number of six digits in the fractional part of the sexagesimal numbers yields an underestimate of the solution.  Using $x_{-1}=1$ and $x_0=2$, iterated linear interpolation (the point 2 will remain fixed) with 6 fractional sexagesimal digits and 14 iterations give
\(1^\circ\;22'\;7''\;42'''\;33^{IV}\;4^V\; 38^{VI}\) and not the sexagesimal digit \(40^{VI}\). The number of iterations corresponds to implementation with 64 bits floating point arithmetic. Based on this slow convergence Maruszewski claims that it is unlikely that Fibonacci used iterated linear interpolation and suggests a fixed point iteration \cite{Maruszewski2009}.

Computer simulation with linear interpolation and sexagesimal computation using a fixed number of digits in the fractional part of the sexagesimal numbers in the formula for linear interpolation does not support the use of linear interpolation where the iterates brackets the solution.

However, working with a dynamic number of digits in the fractional part Glushkov \cite{Glushkov1976} shows using the secant method:
\begin{itemize}
\item Step 1 $x_{-1}=1 < x_0=2$. Computing the new point $x_1$ with one fractional sexagesimal digit $x_1=1^\circ\;18' (<x_*)$
\item Step 2 $x_{1}< x_0$. Computing the new point $x_2$ with one fractional sexagesimal digit $x_2=1^\circ\;20' (<x_*)$
\item Step 3 $x_{1}< x_2$. Computing the new point $x_3$ with one more fractional sexagesimal digit $x_3=1^\circ\;22'\;58'' (>x_*)$ ($x_{1}$ and  $x_2$ is not bracketing the solution)
\item Step 4 $x_{2}< x_3$. Computing the new point $x_4$ with two fractional sexagesimal digits $x_4=1^\circ\;21'\;13'' (<x_*)$
\item Step 5 $x_{4}< x_3$. Computing the new point $x_5$ with one more fractional sexagesimal digit $x_5=1^\circ\;22'\;06''\;32''' (<x_*)$
    \end{itemize}
After 18 iterations Glushkov\footnote{With full accuracy in $p(x_i)$ there is a minor error in step 3. Correct values are $x_3=1^\circ\;22'\;10$, $x_4=1^\circ\;22'\;7''$, $x_5=1^\circ\;22'\;7''\;42'''$, $x_6=1^\circ\;22'\;7''\;42'''\;33^{IV}$, and $x_7=1^\circ\;22'\;7''\;42'''\;33^{IV}\;4^{V}\;38^{VI} $. If two truncated iterates are equal, an additional digit is added.}
 gets Fibonacci's result (\ref{eq:Fib_number}). In the secant method we do not restrict the points to bracket the solution. The convergence rate of the secant method is around 1.62 compared to 2 for the Newton-Raphson method. The rounding performed in Glushkov's implementation of the secant method impedes the rate of convergence compared to Gram's implementation of the Newton-Raphson method.

The father of the History of Science, George Sarton, claims that Fibonacci must have used Regula Falsi (iterated linear interpolation) \cite[p.611]{Sarton1931}.
Fibonacci demonstrates in {\em Liber abaci} the use of linear interpolation on linear problems and Brown \cite{Brown2008} concludes in 2008 that Fibonacci is using iterated linear interpolation also for the nonlinear cubic equation.
%\section{Linear Interpolation}\label{sec:interpolation}
\section{Geometric solution}\label{sec:geometric}
Let $y= \frac{2\sqrt{10}}{x}$ and consider
\[x^2(x^2+(y-\sqrt{10})^2 - 4) = (x^3+2x^2+10x-20)(x-2).\]
If $x$ is a solution of (\ref{eq:Fib}) then $(x,y)$ is the intersection of the circle
$x^2+(y-\sqrt{10})^2 = 4$   and hyperbola $xy=2\sqrt{10}$. This is shown in Figure \ref{fig:Geometric}.

Omar Khayyam (1048--1131) is the first to give a general theory of cubic equations and the first to geometrically solve every type of cubic equation for positive roots.
Geometric solution of (\ref{eq:Fib_depressed}) is the intersection of a parabola and a semicircle shown in Figure \ref{fig:Geometric_depressed}.
For the geometric construction
\cite{Woepcke1851,Kasir1931,Amir-Moez1961,Amir-Moez1962,Lumpkin1978,Rizvi1985,Linden2012,Kent2016,Siadat2021}.

\begin{figure}[ht]
  \includegraphics[scale=.8]{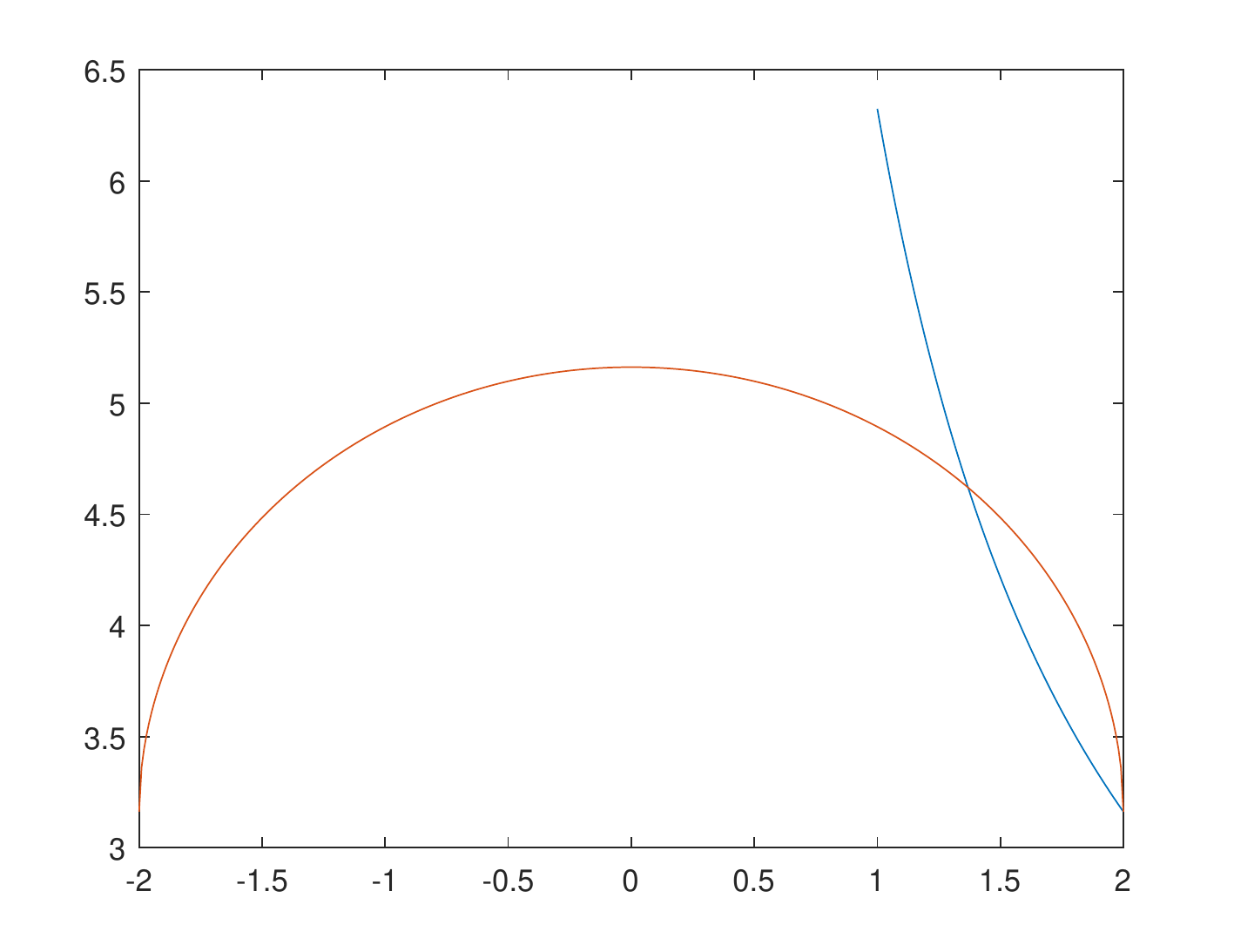}
\caption{Hyperbola $xy=bc$ and semicircle $(x-\frac{c-a}{2})^2+(y-b)^2=\left(\frac{a+c}{2}\right)^2$ where $a=2$, $b=\sqrt{10}$ and $c=2$. The intersection is the positive solution of $x^3+ax^2+b^2x=b^2c$ or $x^3+2x^2+10x=20.$}
\label{fig:Geometric}
\end{figure}
\begin{figure}[ht]
\includegraphics[scale=.8]{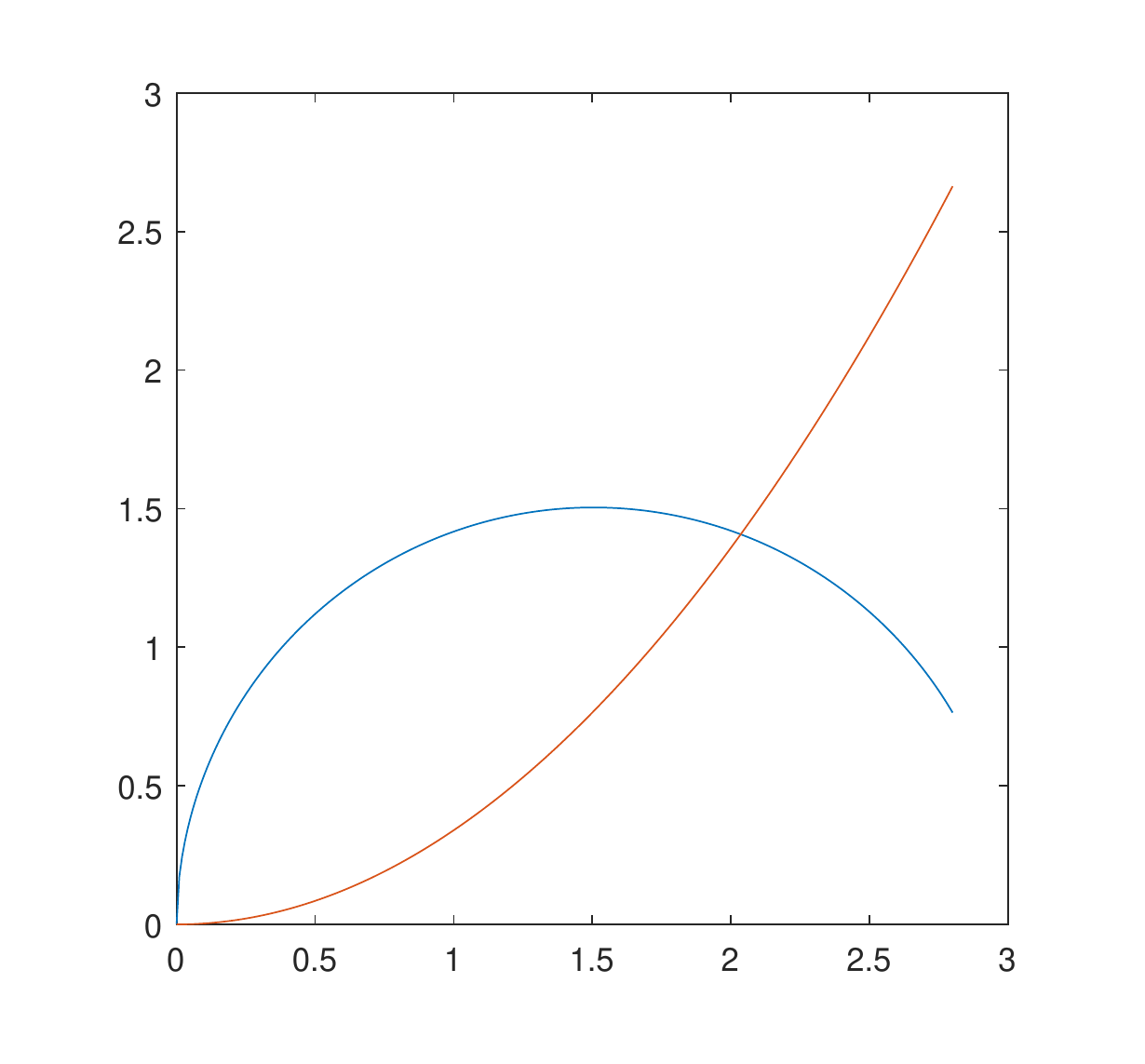}
\caption{Parabola $y=\frac{1}{b}x^2$ and semicircle with diameter $c$, $(x-\frac{c}{2})^2+y^2=\left(\frac{c}{2}\right)^2$ where $b=\sqrt{\frac{26}{3}}$ and $c= \frac{352}{117}.$ The intersection is the positive solution of $x^3+b^2 x=b^2c$ or $x^3+\frac{26}{3}x=\frac{704}{27}$}
\label{fig:Geometric_depressed}
\end{figure}

The work of Omar Khayyam was definitely known to Fibonacci and the geometric solution will easily yield the first fractional sexagesimal digit. This would shorten the computation using the digit--by--digit approach.

\section{Concluding remarks}

We claim that the algorithm used by Fibonacci most likely is a digit--by--digit method. The method was well known at the time of Fibonacci, however we do not have any manuscripts giving details of an implementation of the algorithm. We have also shown that it is possible to get Fibonacci's result using what today is called the Newton--Raphson method with suitable truncation or rounding. The binomial theorem was known at the time of Fibonacci and also the linear approximation/correction. The Newton-Raphson method was described in the late 17th century (1685 and 1690), but the linear approximation was used long before. We have also shown that Fibonacci's result can be obtained from iterated linear interpolation. Linear interpolation is described by Fibonacci in {\em Liber abaci} and thus well known. The single argument that Fibonacci used digit--by--digit is that all arguments are indirect. However, it is the method with the least number of arithmetic operations.

Plausible other methods exist. Maruszewski suggests a fixed--point iteration based on Heron's method \cite{Maruszewski2009}.

In reproducing the method we also need to know the rounding used by the mathematicians.
Reproducing results published in scientific papers or even textbooks may be a challenging exercise for college students and teachers. In this paper we find four cited papers addressing Fibonacci's solution \cite{Brown2008,Maruszewski2009} and Khayyam's geometric solution of cubic equations \cite{Kent2016,Lumpkin1978} suitable for college students.

\bibliographystyle{plain}
%\bibliography{Test_bibliography}
\bibliography{Fibonacci,../Test_bibliography_test}

\begin{thebibliography}{10}

\bibitem{Amir-Moez1962}
A.~R. Amir-Moèz.
\newblock Khayyam's solution of cubic equations.
\newblock {\em Mathematics Magazine}, 35(5):269--271, 1962.

\bibitem{Amir-Moez1961}
A.~R. Amir-Moèz.
\newblock A paper of {O}mar {K}hayyam.
\newblock {\em Scripta Mathematica}, 26(5):323--337, 1962.

\bibitem{Amir-Moez1994}
A~R Amir-Moéz.
\newblock {K}hayyam, al-{B}iruni, {G}auss, {A}rchimedes, and quartic equation.
\newblock {\em The Texas journal of science}, 46:241--257, 1994.

\bibitem{Flos1854}
Baldassarre Boncompagni.
\newblock {\em {Tre scritti inediti di {L}eonardo {P}isano}. Transcribed and
  translated by {B}aldassarre {B}oncompagni}.
\newblock Firenze, 1854.

\bibitem{Brown2008}
Ezra Brown and Jason~C. Brunson.
\newblock Fibonacci's forgotten number.
\newblock {\em College Mathematics Journal}, 32(2):112--120, 2008.

\bibitem{Calinger1999}
Ronald Calinger.
\newblock {\em A Contextual History of Mathematics: To Euler}.
\newblock Prentice Hall, 1999.

\bibitem{Cantor1892}
Moritz Cantor.
\newblock {\em Vorlesungen über geschichte der mathematik, Zweiter Band}.
\newblock Leibzig, 1892.

\bibitem{Cantor1900}
Moritz Cantor.
\newblock {\em Vorlesungen über geschichte der mathematik, Zweiter Band, 2nd
  edition}.
\newblock Leibzig, 1900.

\bibitem{Cassina1924}
Ugo Cassina.
\newblock Risoluzione graduale dell'equazione di {L}eonardo {P}isano.
\newblock {\em Atti della Accademia delle Scienze di Torino}, 59:14--29, 1924.

\bibitem{Coolidge1949}
J.~L. Coolidge.
\newblock The story of the binomial theorem.
\newblock {\em The American Mathematical Monthly}, 56(3):147--157, 1949.

\bibitem{Dattolini2021}
Giuseppe Dattoli, Silvia Licciardi, and Marcello Artioli.
\newblock {\em Vedic Mathematics: A Mathematical Tale From The Ancient Veda To
  Modern Times}.
\newblock World Scientific Publishing, Singapore, 2021.

\bibitem{Drozdyuk2010}
Andriy Drozdyuk and Denys Drozdyuk.
\newblock {\em Fibonacci, his numbers and his rabbits}.
\newblock Choven Publishing, Toronto, 2010.

\bibitem{Genocchi1855}
Angelo Genocchi.
\newblock Sopra tre scritti inediti di {L}eonardo {P}isano pubblicati da {B.}
  {B}oncompagni.
\newblock {\em Annali di scienze mathematiche e fisiche}, 6:161--185, 218--251,
  273--320, 345--362, 1855.

\bibitem{Glushkov1976}
Stanislaw Glushkov.
\newblock On approximation methods of {L}eonardo {F}ibonacci.
\newblock {\em Historia Mathematica}, 3(3):291--296, 1976.

\bibitem{Gram1893}
J{\o}rgen~Pedersen Gram.
\newblock Essai sur la restitution du calcul de {L}\'{e}onard de {P}ise sur
  l'\'{e}quation $x^3+2x^2+10x=20$.
\newblock {\em Kongelige Danske Videnskabernes Selskabs Forhandlinger}, pages
  18--28, 1893.

\bibitem{Hankel1874}
Hermann Hankel.
\newblock {\em Geschichte der Mathematik in Alterthum und Mittelalter}.
\newblock Leipzig, 1874.

\bibitem{Harriot1631}
Thomas Harriot.
\newblock {\em {Artis Analyticae Praxis}}.
\newblock London, 1631.

\bibitem{Kasir1931}
Daoud~Suleiman Kasir.
\newblock {\em The Algebra of {O}mar {K}hayyam}.
\newblock PhD thesis, Columbia University, New York City, 1931.

\bibitem{Kent2016}
Deborah~A. Kent and David~J. Muraki.
\newblock A geometric solution of a cubic by {O}mar {K}hayyam… in which
  colored diagrams are used instead of letters for the greater ease of
  learners.
\newblock {\em The American Mathematical Monthly}, 123(2):149--160, 2016.

\bibitem{Woepcke1851}
Omar Khayyam.
\newblock {\em {L'Algèbre d'Omar Alkhayyâmî}. Translated and commented by F.
  Woepcke}.
\newblock Springer, 1851.

\bibitem{Linden2012}
Sebastian Linden.
\newblock {\em {Die Algebra des Omar Chayyam}. 2nd edition}.
\newblock Springer Spektrum, 2017.

\bibitem{Lumpkin1978}
Beatrice Lumpkin.
\newblock A mathematics club project from {O}mar {K}hayyam.
\newblock {\em The Mathematics Teacher}, 71(9):740--744, 1978.

\bibitem{Maruszewski2009}
Richard Maruszewski.
\newblock Fibonacci's forgotten number revisited.
\newblock {\em College Mathematics Journal}, 40(4):248--251, 2009.

\bibitem{Oughtred1647}
William Oughtred.
\newblock {\em {The key of the mathematics new forged and filed}}.
\newblock London, 1647.

\bibitem{Papakonstantinou2013}
Joanna~M. Papakonstantinou and Richard~A. Tapia.
\newblock Origin and evolution of the secant method in one dimension.
\newblock {\em The American Mathematical Monthly}, 120(6):500--518, 2013.

\bibitem{Flos1983}
Ettore Picutti.
\newblock Il {Fl}os di {Leonardo Pisano} dal codice {E.75.P.} sup. della
  {Biblioteca Ambrosiana di Milano}.
\newblock {\em Physis: Rivista Internazionale di Storia della Scienza},
  25(2):293--387, 1983.

\bibitem{Posamentier2007}
Alfred~S. Posamentier and Ingmar Lehmann.
\newblock {\em The Fabulous Fibonacci Numbers}.
\newblock Prometheus Books, New York, 2007.

\bibitem{Rizvi1985}
Vaqar~Ahmed Rizvi and Y.A Rizvi.
\newblock {'U}mar {K}hayyām as a geometrician—a survey.
\newblock {\em Islamic Studies}, 24(2):193--204, 1985.

\bibitem{Roberts2020}
Alexandre~M. Roberts.
\newblock Mathematical philology in the treatise on double false position in an
  {A}rabic manuscript at {C}olumbia university.
\newblock {\em Philological Encounters}, 5(3-4):308 -- 352, 2020.

\bibitem{Sarton1931}
George Sarton.
\newblock {\em Introduction to the History of Science: Volume 2}.
\newblock Number Part II in Carnegie Institution of Washington. Publication. no
  376. Williams and Wilkins Company, 1931.

\bibitem{Schwartz2004}
Randy~K. Schwartz.
\newblock Issues in the origin and development of {H}isab al-{K}hata'ayn
  ({C}alculation by double false position).
\newblock Eighth North African Meeting on the History of Arab Mathematics.
  Radès, Tunisia, 2004.

\bibitem{Siadat2021}
M.~Vali Siadat and Alana Tholen.
\newblock {O}mar {K}hayyam: {G}eometric algebra and cubic equations.
\newblock {\em Math Horizons}, 28(1):12--15, 2021.

\bibitem{Smith1925}
David~Eugene Smith.
\newblock {\em History of mathematics. Volume II. Special topics of elementary
  mathematics}.
\newblock Ginn and Company, Boston, 1925.

\bibitem{Steihaug2021}
Trond Steihaug.
\newblock Computational science in the 17th century. numerical solution of
  algebraic equations: Digit–by–digit computation.
\newblock In Mehiddin Al-Baali, Anton Purnama, and Lucio Grandinetti, editors,
  {\em Numerical Analysis and Optimization, NAO-V, Muscat, Oman, January 2020},
  chapter~12, pages 249--269. Springer Nature, 2021.

\bibitem{Stevin1585}
Simon Stevin.
\newblock {\em L'Arithmetique, Contenant les computations des nombres
  Arithmetiques ou vulgaires : Au{\ss}i l'Algebre ... Ensemble les quatre
  premiers livres d'Algebre de Diophante d'Alexandrie ... traduicts en Francois
  ; Encore un livre particulier de la Pratique d'Arithmetique (etc.)}.
\newblock Christophle Plantin, Leiden, The Netherlands, 1585.

\bibitem{Terquem1856}
Olry Terquem.
\newblock Sur {L}éonard {B}onacci de {P}ise et sur trois écrits de cet
  autheur publiés par {B}althasar {B}oncompagni.
\newblock {\em Annali di scienei mathematiche e fisiche}, 7:106--147, 1856.

\bibitem{Tropfke1980}
Johannes Tropfke.
\newblock {\em Geschichte der elementarmathematik, Band 1. Arithmetik und
  algebra}.
\newblock Walther de Gruytter, Berlin, 4th edition, 1980.

\bibitem{Vetter1928}
Quido Vetter.
\newblock Nota alla risoluzione dell'equazione cibica di {L}eonardo {P}isano.
\newblock {\em Atti della Accademia delle Scienze di Torino}, 63:296--299,
  1928.

\bibitem{Viete1600}
Fran\c{c}ois Vi\`{e}te.
\newblock {\em {De numerosa potestatum ad exegesim resolutione}}.
\newblock Paris, 1600.

\bibitem{WhitesideI}
Derek~Thomas Whiteside.
\newblock {\em The Mathematical Papers of {I}saac {N}ewton 1664-1666, {V}olume
  I}.
\newblock Cambridge University Press, Cambridge, 1967.

\bibitem{Woepcke1854}
Franz Woepcke.
\newblock Sur un essai de déterminer la nature de la racine d’une équation
  du troisième degré, contenu dans un ouvrage de {L}éonard de {P}ise
  découvert par {M}. le prince {B}althasar {B}oncompagni.
\newblock {\em Journal de mathématiques pures et appliquées}, 1re série,
  tome 19:401--406, 1854.

\bibitem{Yadegari1980}
Mohammad Yadegari.
\newblock The binomial theorem: {A} widespread concept in medieval {I}slamic
  mathematics.
\newblock {\em Historia Mathematica}, 7(4):401--406, 1980.

\bibitem{Ypma1995}
Tjalling~J. Ypma.
\newblock Historical development of the {N}ewton-{R}aphson method.
\newblock {\em SIAM Review}, 37(4):531--551, 1995.

\bibitem{Zeuthen1893}
Hieronymus~Georg Zeuthen.
\newblock Sur la r\'{e}solution numerique d'un\'{e}equation du $3^e$ degr\'{e}
  par {L}\'{e}onard de {P}ise.
\newblock {\em Kongelige Danske Videnskabernes Selskabs Forhandlinger}, pages
  1--17, 1893.

\end{thebibliography}
\end{document}